\documentclass[10pt]{amsart}

\usepackage{amscd}
\usepackage{amssymb} 

\def\BbbC{\Bbb C}

\def\build#1_#2^#3{\mathrel{\mathop{\kern 0pt#1}\limits_{#2}^{#3}}}

\def\build#1_#2^#3{\mathrel{\mathop{\kern 0pt#1}\limits_{#2}^{#3}}}

\def\N{\nabla}
\newcommand{\W}{{\mathcal W}}
\newcommand{\R}{{\mathcal R}}
\let\dd=\displaystyle

\newtheorem{theoreme}{Theorem}[section]
\newtheorem{definition}[theoreme]{Definition}
\newtheorem{proposition}[theoreme]{Proposition}

\newtheorem{lemme}[theoreme]{Lemma}
\newtheorem{corollaire}[theoreme]{Corollary}

\let\dd=\displaystyle

\def\hfleche#1#2{\smash{\mathop{\vbox{\hbox to
        12mm{{\rightarrowfill}}}}\limits^{#1}_{#2}}}

\title{\sc Linearizable 3-webs and the Gronwall conjecture}

\author{Joseph GRIFONE} 
\address{J.~GRIFONE: Laboratoire Emile Picard, U.M.R.  C.N.R.S. 5580,
  D\'epartement de Math\'ematiques, Universit\'e Paul Sabatier, 118,
  Route de Narbonne, 31062 Toulouse Cedex, France}
\email{grifone@picard.ups-tlse.fr}

\author{Zoltan MUZSNAY} 
\address{Z.~Muzsnay: Department of Mathematics, University of
  Debrecen, Debrecen, H-4032 PBox 12, Hungary}
\email{muzsnay@math.klte.hu}

\author{Jihad SAAB}
\address{J.~SAAB: Facult\'e des sciences et de Genie, Universit\'e
  Saint Esprit de Kaslik, B.P. 446 Jounieh, Lebanon}
\email{jihadsaab@yahoo.fr}


\thanks{The authors thank P.\,T.\,Nagy for the many discussions which
  made progress possible in this work. They also thank A. \,Henaut for
  his encouragements and valuable consultations.}

\begin{document}

\maketitle \vskip 3 cm

\begin{abstract}  
  In the article \cite{GMS} published in 2001 in the journal
  "Nonlinear Analysis", we studied the linearizability problem for
  $3$-webs on a $2$-dimensional manifold.  Four years after the
  publication of our article Goldberg and Lychagin \cite{GL_cras}
  obtained similar results by a different method and criticized our
  article by qualifying the proofs incomplete. However, they obtained
  false result on the linearizability of a certain web.  We present
  here the complete version of \cite{GMS} with computations and
  explicit formulas, because we deem that the opinion of Goldberg and
  Lychagin in \cite{GL_cras} concerning our work is unjustified.
\end{abstract} 

\medskip
\medskip

\section{Introduction}

In the article \cite{GMS} published in 2001 in the journal "Nonlinear
Analysis", we studied the linearizability problem for $3$-webs on a
$2$-dimensional manifold.  Using the integrability theory of
over-determined partial differential systems, we computed the
obstructions to linearizability and we produced an effective method to
test the linearizability of 3-webs in the (real or complex) plane.  We
showed that, in the non-parallelizable case, there exists an algebraic
submanifold ${\mathcal A}$ of the space of vector valued symmetric
tensors ($S^2T^*\otimes T$), which can be expressed in terms of the
curvature of the Chern connection and its covariant derivatives up to
order $6$, such that the affine deformation tensor is a section of
$S^2T^*\otimes T$ with values in ${\mathcal A}$.  In particular, we
proved that a web is linearizable if and only if ${\mathcal A}\neq
\emptyset$, and there exists at most $15$ projectively nonequivalent
linearizations of a nonparallelizable 3-web.  In order to give a
coordinate free and intrinsic presentation of the results we used
tensors and covariant derivatives to find the obstructions to the
linearization.

Recently Goldberg and Lychagin \cite{GL_cras} obtained similar results
by a different method. They criticized our article by qualifying the
proofs incomplete, without giving any justification or reason for
their claim. They claim that "...the main and only example of a
linearizable (in their approach) 3-web ... is not linearizable at
all..."  To prove their statement they apply their theory to this
particular web and find that the corresponding algebraic submanifold
is empty.  However, in the article \cite{MZ} which appears in arXiv
with this present paper, Z.~Muzsnay shows by producing an explicit
linearization, that in accordance with the claim of \cite{GMS}, this
web is linearizable. This proves that something is wrong in their
work: either the proofs of \cite{GL} and \cite{GL_cras} are not
correct, or some of their calculations are false.
  
We are putting on the arXiv preprint archive a detailed version of the
article \cite{GMS} with computations and explicit formulas, because we
deem that the opinion of Goldberg and Lychagin in \cite{GL_cras}
concerning our work is unjustified.

\bigskip

\section{Introduction to the linearizability problem of 3-webs}

Let $M$ be a two-dimensional real or complex differentiable manifold.
A $3$-web is given in an open domain $D$ of $M$ by three foliations of
smooth curves in general position.  Two webs ${\mathcal W}$ and
$\tilde{{\mathcal W}}$ are locally equivalent at $p\in M$, if there
exists a local diffeomorphism on a neighborhood of $p$ which exchanges
them.

A $3$-web is called {\it linear} (resp. {\it parallel}) if it is given
by 3 foliations of straight lines (resp. of parallel lines). A
$3$-web which equivalent to a linear (resp. parallel) web is called
linearizable (resp. parallelizable).

A linear connection, called Chern connection and denoted by $\nabla$,
can be associated to a $3$-web. $\nabla$ preserves the web, i.e.  the
leaves are auto-parallel curves. It is not difficult to see that a
$3$-web is parallelizable if and only if the curvature of the Chern
connection vanishes.  The Graf-Sauer Theorem (\cite{BB}, page 24)
gives an elegant characterization of such webs: a linear web is
parallelizable if and only if, its leaves are tangent lines to a curve
of degree 3. 

\smallskip

The problem to give linearizability criterion is a very natural one.
Such criterion is important in nomography (cf.  \cite{Gro}):
determining whether some nomogram can be reduced to an alignment chart
is equivalent to the problem of determining whether a web is
linearizable.  The most significant works on this subject are due to
Bol (\cite{Bol1}, \cite{Bol2}).  In \cite {Bol1} he suggested how to
find a criterion of linearizability, although he is unable to carry
out the computation, which really need the use of computer.  He shows
that the number of projectively different linear $3$-webs in the plane
to which a non-hexagonal $3$-web is equivalent is finite and less that
17.  Bol's proof consists in to associate to a real $3$-web two
complex vector fields which play an essential role, so his proof
cannot be translated in the complex case. In our computation the web
can be real as well as complex.

The formulation of the linearizability problem in terms of Chern connection
was suggested by Akivis in a lecture given in Moscow in 1973.  Following
Akivis idea Goldberg in \cite{Gol2} found all affine connection $\Gamma^*$
relative to which the web leaves are geodesic lines and distinguished a
linearizable 3-webs by claiming that the connection $\Gamma^*$ is flat.  In
this paper we are using this approach to solve the problem. 

Denoting by $T$ and $T^*$ the tangent and the cotangent bundle of $M$, a
section $L$ of the bundle $S^2T^*\otimes T$ on $M$ is called {\it
pre-linearization}, if the connection $\nabla^L$ defined by
\begin{displaymath}
  \nabla^L_XY = \nabla_XY + L(X,Y)
\end{displaymath}
preserves the web, that is the three families of leaves are
auto-parallels curves with respect to $\N^L$.  A pre-linearization $L$
is called {\it linearization} if the connection $\N^L$ is flat i.e.
the curvature of the connection $\N^L$ given by equation (\ref{P1})
vanishes.  This equation gives us a first order partial differential
system on $L$.  Two linearizations $L$ and $L'$ are projectively
equivalent if the connections $\N^{L}$ and $\N^{L'}$ are projectively
related. The equivalences classes are called classes of
linearizations. They are in one-to-one correspondence with the {\it
  bases of the linearization} which is a simple projective invariant,
noted by $s$. The linearizability condition can be reformulate with
this object by a second order partial differential system.  We show
that the system is of finite type, and the obstruction to the
linearizability can be expressed in terms of polynomials of $s$, whose
coefficients depends only on the curvature tensor of the Chern
connection.  Our main result is the following:

\smallskip
\begin{theoreme}
  Let $\W$ be an analytical $3$-web on a 2-dimensional real or complex
  manifold $M$, whose Chern curvature does not vanish at $p\in M$.  Then,
  there exists an algebraic sub-manifold ${\mathcal A}$ of $E$ over a
  neighborhood of $p$, expressed in terms of the curvature of the Chern
  connection and its covariant derivatives up to order $6$, so that the
  linearizations of $\W$ are sections of $E$ with values in ${\mathcal A}$.
  In particular:
  \begin{enumerate}
  \item The web is linearizable if and only if ${\mathcal A}\neq
    \emptyset$;
  \item There exists at most $15$ classes of linearizations.
  \end{enumerate}
  The explicit expression of the polynomials and its coefficients which
  define $\mathcal{A}$ can be found in Chapter 6 and 7.
\end{theoreme}

\bigskip

\section{Notations and definitions}

\noindent
Let ${\mathcal W}$ be a differential 3-web on a manifold $M$ given by a triplet
of mutually transversal foliations $\{{\mathcal F}_1,{\mathcal
  F}_2,{\mathcal F}_3\}$.
From the definitions it follows that $M$ is even dimensional and that the
dimension of the tangent distributions of the foliations ${\mathcal
F}_1$, ${\mathcal F}_2$, ${\mathcal F}_3$ is the half of the dimension of $M$.
The foliations $\{{\mathcal F}_1,{\mathcal F}_2,{\mathcal F}_3\}$ are
called horizontal, vertical and transversal and their tangent space are
denoted by $T^h$, $T^v$ and $T^t$.

The following theorem proved by Nagy \cite{Nag} gives an elegant
infinitesimal characterization of $3-$webs and their Chern connection.

\smallskip

\noindent
\textbf{Theorem.}  \textit{A $3-$web is equivalent to a pair $\{h,j\}$
  of (1,1)-tensor fields on the manifold, satisfying the following
  conditions:
  \begin{enumerate}
  \item $h^2 = h$, $j^2 = id$,
  \item $jh = vj$, where $v=id-h$,
  \item ${\rm Ker}\,h$, ${\rm Im}\, h$ and ${\rm Ker}(h + id)$ are
    integrable distributions.
  \end{enumerate}
  For any $3-$web, there exists a unique linear connection $\N$ on $M$
  which satisfies
  \begin{enumerate}
  \item $\nabla h = 0$,
  \item $\nabla j = 0$,
  \item $T(hX,vY) = 0$, \ for every $X,Y\in TM$, $T$ being the
    torsion tensor of $\N$.
  \end{enumerate} 
  $\nabla$ is called {\rm Chern connection}.}

\medskip

\noindent
In the sequel, we suppose that the dimension of $M$ is two.
\begin{definition}
  Let $\W$ be a $3$-web and $\N$ its Chern connection. A symmetrical
  (1,2)-tensor field $L$ is called pre-linearization if the
  connection
  \begin{displaymath}
    \N^L_XY = \N_XY + L(X,Y)
  \end{displaymath}
  preserves the web, that is the leaves are auto-parallel curves with
  respect to $\N^L$.  A pre-linearization is a linearization if the
  connection $\N^L$ is flat i.e. its curvature vanishes.  Two
  pre-linearizations $L$ and $L'$ are projectively equivalent if the
  connections $\N^L$ and $\N^{L'}$ are projectively related, that is
  there exists $\omega\in \Lambda^1(M)$, such that
  $$\N^L_XY = \N^{L'}_XY + \omega(X) Y + \omega(Y) X$$
\end{definition}

\begin{proposition}
  \label{def_linearisation}
  A tensor field $L$ in $S^2T^* \otimes T$ is a {\rm linearization} if
  and only if
  \begin{alignat*}{1}
    1. & \ vL(hX,hY) = 0,
    \\
    2. & \ hL(vX,vY) = 0,
    \\
    3. & \ L(hX,hY)+jL(jhX,jhY)-hL(jhX,hY)
    \\
    & \qquad \qquad \qquad - hL(hX,jhY) -jvL(jhX,hY)-jvL(hX,jhY) = 0,
    \\
    4. & \ \N_XL (Y,Z) - \N_YL (X,Z) + L(X,L(Y,Z)) - L(Y,L(X,Z)) + R
    (X,Y)Z = 0.
  \end{alignat*}
  holds, for any $X,Y,Z \in T$, where $R$ denotes the curvature of the
  Chern connection.
\end{proposition}
The proof is a straightforward verification. Properties 1), 2) and 3) means
that $L$ is a pre-linearization and follows from the fact that $\N^L$
preserves the web, while properties 4) expresses, that the curvature of
$\N^L$ vanishes.

\smallskip

\begin{definition}
  \label{base}
  Let $M$ be a $2-$dimensional manifold, $\W$ a web on $M$ and
  $\{e_1,e_2\}$ a frame at $p\in M$ adapted to the web, i.e. $e_1\in
  T^h_p,\, e_2=je_1 \in T_p^v$. Let $L$ be a pre-linearization at $p$,
  whose components are $L_{ij}^k$, that is: $L(e_i,e_j) =
  L^k_{ij}e_k$, and let us set the tensor-field $s$ represented by
  the components $2L^1_{12} - L^2_{22}$.  The tensor $s$ will be
  called the base of $L$.
\end{definition}  

The following proposition is elementary, but it is the key for the
proof of our main theorem.

\begin{proposition}
  \label{ss'}
  Two pre-linearizations $L$ and $L'$ are projectively equivalent if
  and only if they have the same base, i.e.  $s = s'$.
\end{proposition}

Indeed, if $L$ and $L'$ are two projectively equivalent
pre-linearizations, then there exists $\omega \in T^*$ such that $L'
= L + \omega \odot id$, i.e. in the frame $\{e_1,e_2\}$ :
\begin{displaymath}
      {L'}^1_{11}=L^1_{11}+2\omega_1,
      \qquad
      {L'}^2_{22}=L^2_{22}+2\omega_2,
      \qquad
      {L'}^1_{12}=L^1_{12}+\omega_2
\end{displaymath} 
where $\omega_1$ and $\omega_2$ are the components of $\omega.$ This
system is consistent if and only if ${L'}^1_{12}-L^1_{12}= \frac{1}{2}
({L'}^2_{22}-L^2_{22})$, \ i.e.  $s=s'$.

\bigskip

\section{The linearization operator}

Let $M$ be a 2-dimensional real or complex manifold and $\W$ a 3-web
on $M$.  $\Lambda^kT^*$ and $S^kT^*$ are the bundles of the
$k$-skew-symmetric and symmetric forms. If $B \to M$ is a vector
bundle on $M$, then ${\mathcal S}ec (B)$ will denote the sheaf of the
sections of $B$ and $J_k(B)$ the vector bundle of $k-$jets of the
sections of $B$.

\medskip

In the sequel $E$ will denotes the bundle of the pre-linearizations
and $F := \Lambda^2T^*\otimes T$.  In order to study the
linearizability of $\W$, we will consider the differential operator
$P_1: E \to F$ and study the integrability of the differential system
\begin{math}
  P_1(L) = 0,
\end{math}
where
\begin{equation}
  \label{P1}
    \begin{aligned}
      \bigl(P_1(L)\bigl)(X,Y,Z) \,= \, & (\N_XL) (Y,Z) - (\N_YL) (X,Z) +
      \\
      & + L(X,L(Y,Z)) - L(Y,L(X,Z)) + R (X,Y) Z
  \end{aligned}
\end{equation}
for every $X, Y, Z \in T$.

\medskip

We will use the theory of the formal integrability of Spencer
(\cite{BCG}, \cite{GM}). The notations are those of \cite{GM}, where
is given also an accessible introduction to this theory. In
particular, if $P$ is a quasi-linear operator of order $k$ and $p \in M$,
then $R_{k, p}$ is the bundle of the formal solutions of order $k$ at
$p$, $\sigma_{k+\ell}(P)$ or simply $\sigma_{k+\ell}$ is the symbol of
the $\ell$-th order prolongation $p_{\ell}(P)$ of $P$.  We also denote
$g_{k+\ell} = {\rm Ker} \sigma_{k+\ell}$ and $K = {\rm Coker}\,
\sigma_{k+1}$.

\medskip

Let $L \in E$ a pre-linearization.  We introduce the tensors
\begin{displaymath}
  x, \, y , \, z : T^h  \otimes T^h \to T^h
\end{displaymath}
defined by
\begin{equation}
  \left\{
    \begin{aligned}
      x \,(hX,\,hY) & = L\,(hX,\,hY)
      \\
      y \,(hX,\,hY) & = j L\,(jhX,\,jhY)
      \\
      z \,(hX,\,hY) & = h L\,(hX\,,jhY)
    \end{aligned}
  \right.
\end{equation} 
One denotes $x^2$ the (1,3) tensor defined by $x^2\,(hX,hY,hZ) = x\,
(x\, (hX,\,hY),\,hZ)$. Similarly, we define the product $xy$, $ x^3$
(which is a (1,4) tensor field), etc...

The space of pre-linearizations, $E$ is a 3-dimensional vector bundle over
$M$,  and $x$, $y$, $z$ can be used to parameterize it. However, taking
into account some symmetries of the problem and the Proposition \ref{ss'},
it is better to introduce the tensors $s,t : T^h\oplus T^h \to T^h$ defined
by
\begin{equation}
  \begin{aligned}
    s &= 2z - y
    \\
    t &= \tfrac {1}{2}\,(x + y - 2z)
  \end{aligned}
\end{equation}
and parameterize $E$ by $s$, $t$, $z$ where $s$ is the base  of the web (see
Definition \ref{base}).

\medskip

In order to simplify the notation, we denote by $C_1$ and $C_2$ the
tensor fields $\big(\bigotimes ^{p +1}{T^h}^*\big)\otimes T^h$ defined
by
\begin{equation}\label {C_1C_2}
  \begin{array}{l}
    C_1\,(hX,hX_1,...,hX_p) = (\N_{hX}C)\,(hX_1,...,hX_p)\\
    C_2\,(hX,hX_1,...,hX_p) = (\N_{jhX}C)\,(hX_1,...,hX_p)
  \end{array}
\end{equation}
where $C$ is a tensor field in $\big(\bigotimes ^p {T^h}^*\big)\otimes
T^h$.  By recursion, we introduce the successive covariant derivatives
with the convention that $C_{i_1i_2} := (C_{i_2})_{i_1}$. Thus,
$x_{i_1,...,i_p}$ is the $(1,p+2)$ tensor defined in an adapted frame
by
\begin{displaymath}
  x _{i_1,...,i_p}\, (\underbrace {e_1,...,e_1}_{p\, \mathrm{times}}, 
  hX,hY) = (\underbrace{\N\N\ \cdots \N}_{p\, \mathrm{times}} x )\, 
  (e_{i_1},...,e_{i_p}\,hX,hY).
\end{displaymath}
We denote ${\mathcal R}$ the tensor ${\mathcal R} : T^h \oplus T^h
\oplus T^h \to T^h $ defined by
\begin{equation}\label{calR}
  {\mathcal R} (hX,\,hY)\,hZ = R(jhX,\,hY)\,hZ
\end{equation}
where $R$ is the curvature of the Chern connection. With the above
notation we have
\begin{displaymath}
  (\nabla_i\nabla_jL^l_{i_1,...,i_m})-(\nabla_j\nabla_iL^l_{i_1,...,i_m})=
  R^l_{ijk}L^k_{i_1,...,i_m}-
  R^k_{iji_1}L^l_{k,...,i_m}-\cdots -R^k_{iji_m}L^l_{i_1,...,k}.
\end{displaymath}
In particular
\begin{equation}
  \label{permutation}
  C_{12} - C_{21} = (p-1)\,{\mathcal R}\,C
\end{equation}
for a tensor field $C \in \big(\bigotimes^p {T^h}^*\big)\otimes T^h$.

\medskip

Using these notations, and resolving two equations in $z_1$ and $t_2$
the system $P_1(L) = 0$ can be write as:
\begin{equation}
  \label{systP1}
  \left\{   
    \begin{aligned}
      t_1 & = st + t^2,
      \\
      t_2 & = \tfrac{1}{3} s_1 - \tfrac{2}{3} s_2 + zt - \tfrac
      {1}{3}{\mathcal R},
      \\
      z_1 & = \tfrac{2}{3} s_1 -\tfrac{1}{3} s_2 + zt + \tfrac{1}{3}
      {\mathcal R},
      \\
      z_2 & = -zs + z^2.
  \end{aligned}
  \right.
\end{equation}
Note that $P_1$ is regular because the symbol and his prolongation are
regular maps.  The system (\ref{systP1}) can be seen as a Frobenius system
on the variables $t$ and $z$, and $s$ being a parameter.  By the formula
(\ref{permutation}), the integrability conditions are
\begin{displaymath}
\left\{
  \begin{aligned}
    z_{12} - z_{21} &= {\mathcal R} z,
    \\
    t_{12} - t_{21} &= {\mathcal R} t,
    \\
    s_{12} - s_{21} &= {\mathcal R} s,
  \end{aligned}
\right.
\end{displaymath}
and thus from (\ref{systP1}) we can arrive to the system
\begin{eqnarray}
  \label{systP2}
  P_2 = \left \{
    \begin{aligned}
      s_{22}=2s_{21}-ss_2+2ss_1+\mathcal{R}s+\mathcal{R}_2,
      \\
      s_{11}=2s_{21}-2ss_2+ss_1+\mathcal{R}s+\mathcal{R}_1.
    \end{aligned}
  \right.
\end{eqnarray}
The operator $P_2 : Sec \,(E_2) \to F_2$ is a quasi-linear second order
differential operator, where $E_2 = {T^h}^*\otimes {T^h}^*\otimes {T^h}$,
and $F_2:=F'\oplus F'$ with $F' := {T^h}^*\otimes {T^h}^*\otimes E_2$.  The
linearizability of the web is equivalent to the integrability of the
operator $P_2$.  In the sequel we will consider this one and examine its
integrability.

\begin{proposition}
  At every $p \in M$ all 2$^{nd}-$order solution at $p$ of $P_2$ can
  be lifted into a 3$^{rh}-$order solution.
\end{proposition} 

Indeed, fixing an adapted base $\{e_1,e_2=je_1\}$, the symbol of $P_2$
is a map $\sigma_2 : S^2T \otimes E_2 \to F_2$ defined by
\begin{displaymath}
  \sigma_2(A) = (A_{22} - 2A_{21},\, A_{11} - 2A_{21}),
\end{displaymath}
where $A_{ij} = A(e_i,e_j)$. So $g_2 := {\rm Ker}\, \sigma_2$ is defined by
the equations $A_{22}-2A_{21} = 0$ and $A_{11}-2A_{21} =0$.  Since these
equations are independent, we have: ${\rm rank} \, \sigma_2 = 2$ and ${\rm
dim}\, g_2 = 1$.  On the other hand, for the first prolongation $\sigma_3:
S^3T^* \otimes E_2 \to T^*\otimes F_2$ we find that $ g_3 = {\rm Ker}\,
\sigma_3$ is defined by the equations
\begin{displaymath}
  B_{k22} - 2B_{k21} = 0, \quad B_{k11} - 2B_{k21} = 0,
\end{displaymath}
$k=1,2$. It is easy to verify that these equations are also independent.
Therefore ${\rm rank} \, \sigma_3 = 4 = {\rm dim} (T^*\otimes F_2)$, and ${\rm
dim} g_3 = 0$, thus $\sigma_3$ is onto i.e. $\textrm{Coker} \, \sigma_3 =
0$.  We have the following exact diagram:
\begin{displaymath}
  \begin{CD}
    && S^{2}T^*\otimes E_2 @>{\sigma _{3}}>> T^*\otimes F @>{ }>> 
    {\dd \frac {T^*\otimes F}{{\rm Im} \sigma _{3}}}=0 
    \\
    && @V{\varepsilon}VV @V{\varepsilon}VV
    \\
    R_{3} @>>> J_{3} (E_2) @>{p_1(P_2)}>> J_1F
    \\
    @V{\overline{\pi}}VV @V{\pi_2}VV @V{\pi _1}VV
    \\
    R_2 @>>> J_2(E_2) @>{p_0(P_2)}>> F
   \end{CD}
 \end{displaymath}
 Consequently $\bar{\pi}_3$ is onto, i.e. every 2$^{nd}-$ order
 solution of $P_2$ can be lifted into a 3$^{rd}-$order solution.  $\Box$

\begin{proposition}
  The operator $P_2$ is not $2-$acyclic, i.e. there is a higher order
  obstruction which arises for the operator $P_2$.
\end{proposition}

Indeed, the sequence
\begin{math}
  0 \longrightarrow g_{{\ell}+1}(P_2) \longrightarrow
  g_{\ell}(P_2)\otimes T^* \stackrel{\delta_{\ell} (P_2)}
  {\longrightarrow} g_{{\ell}-1}(P_2)\otimes
  \Lambda^2T^*\longrightarrow 0
\end{math}
is not exact for all $l \geq 2$, where $\delta_{\ell}$ denotes the
skew-symmetrization in the corresponding variables: for ${\ell = 3}$
we have ${\rm rank}\, \delta_3 = 0 < {\rm dim}
(g_2\otimes{\Lambda}^2T^*)=1$.  $\Box$

\bigskip

\section{The first obstruction}

\noindent
In order to find the higher order obstruction we consider the prolongation
of $P_2$, i.e.  the operator $P_3:=(P_2, \nabla P_2)$, where
\begin{math}
  \nabla P_2: T^* \otimes E_2 \, \longrightarrow \, T^*\otimes F_2
\end{math}
is the covariant derivative of $P_2$ with respect to the Chern
connection.  Explicitly, this system is formed by the system
(\ref{systP2}) and by the following equations:
\begin{eqnarray}
\label{systP3}
\left\{
  \begin{aligned}
    s_{212}=&ss_{21}-\tfrac{1}{3}s_1s_2+\tfrac{4}{3}s_2^2 -
    \tfrac{2}{3}s_1^2+ \tfrac{4}{3}\mathcal{R}s_2 +2s^2s_1
    \\
    & +\mathcal{R}s^2+(2\mathcal{R}_2-\mathcal{R}_1)s-
    \tfrac{2}{3}\mathcal{R}_{21} -\tfrac{1}{3}\mathcal{R}_{12}
    \\
    s_{211}=&-ss_{21}+\tfrac{1}{3}s_1s_2+\tfrac{2}{3}s_2^2
    -\tfrac{4}{3}s_1^2+ (\tfrac{5}{3}\mathcal{R}+2s^2)s_2
    \\
    & -10\mathcal{R}s_1+(\mathcal{R}_2-2\mathcal{R}_1)s
    -\tfrac{1}{3}\mathcal{R}_{21} -\tfrac{2}{3}\mathcal{R}_{12}
    \\
    s_{111}=&-2ss_{21}-\tfrac{4}{3}s_1s_2 + \tfrac{4}{3}s_2^2 -
    \tfrac{5}{3}s_1^2+ (\tfrac{10}{3}\mathcal{R} +2s^2)s_2
    -(\tfrac{5}{3}\mathcal{R}-s^2)s_1
    \\
    & -\mathcal{R}s^2+(2\mathcal{R}_2- 2\mathcal{R}_1)s -
    \tfrac{2}{3}\mathcal{R}_{21} -\tfrac{4}{3} \mathcal{R}_{12} +
    \mathcal{R}_{11}
    \\
    s_{222}=&2ss_{21}+\tfrac{4}{3}s_1s_2+\tfrac{5}{3}s_2^2 -
    \tfrac{4}{3}s_1^2 + (\tfrac{5}{3}\mathcal{R}
    +s^2)s_2-(\tfrac{10}{3}\mathcal{R}-2s^2)s_1
    \\
    &+\mathcal{R}s^2 + (2\mathcal{R}_2- 2\mathcal{R}_1)s
    -\tfrac{4}{3}\mathcal{R}_{21} -\tfrac{2}{3}\mathcal{R}_{12} +
    \mathcal{R}_{22}
  \end{aligned}
\right.
\end{eqnarray}
Since (\ref{systP3}) can be solved with respect to the 3$^{rd}-$order
derivatives, the existence of a 2$^{nd}$-order formal solution implies the
existence of $3^{rd}-$order solutions. 

\smallskip

In the following, we will use the notion of involutivity of a differential
system in the sense used in monograph \cite{BCG} , p. 121 (cf the discussion of
this notion on p. 2) \footnote{Sometimes there is a confusion between
  different terminologies. The involutivity here (and also in the mentioned
  \cite{BCG} and \cite{GM}) means the involutivity of the symbol i.e. that
  the Cartan's test for involutivity holds. It doesn't mean the
  integrability, which is the case in some another terminologies.}. 
We have then
\begin{lemme}
  $P_3$ is involutive.
  Moreover, any $3^{rd}-$order solution of $P_3$ can be lifted into a
  $4^{th}-$order solution if and only if $\varphi =0$, where
  \begin{equation}
    \label{varphi}
    \begin{aligned}
      \varphi(s):=&-24
      \mathcal{R}s_{21}-(24\mathcal{R}s+12\mathcal{R}_1
      -6\mathcal{R}_2)s_1 + (24\mathcal{R}s + 6\mathcal{R}_1 -
      8\mathcal{R}_2)s_2
      \\
      &+3\mathcal{R}s^3 +(-4\mathcal{R}_2-3\mathcal{R}_{22} +
      \mathcal{R}_{21} +2\mathcal{R}_{12}-13 \mathcal{R}^2
      -3\mathcal{R}_{11})s
      \\
      &+2\mathcal{R}_{122}-\mathcal{R}_{221}-\mathcal{R}_{112}
      -5\mathcal{R}\mathcal{R}_1 -2\mathcal{R}_{121} -11\mathcal{R}
      \mathcal{R}_2
    \end{aligned}
  \end{equation}
\end{lemme}

{\it Proof}.  The symbol of $P_3$ is just $\sigma_3(P_3) :
S^3T^*\otimes E_2 \longrightarrow T^* \otimes F_2$, the first
prolongation of the symbol of $P_2$. On the other hand, $\sigma_4(P_3) :
S^4T^* \otimes E_2 \longrightarrow S^2T^*\otimes F_2$ and $g_4 = {\rm
  Ker} \sigma_4$ is defined by the equations
\begin{displaymath}
      D_{ij}^1 := C_{ij22} - 2C_{ij21} = 0,
      \quad
      D_{ij}^2 := C_{ij11} - 2C_{ij21} = 0, \quad i,j =1,2.
\end{displaymath}
There is one relation between these equations: $D_{11}^1 - 2D_{12}^1 -
D_{22}^2 + 2D_{12}^2 = 0$.  Therefore the rank of this system is 5, so if
$K_2$ denotes the cokernel of $\sigma_4$ i.e. $K_2 = (S^2T^*\otimes F_2) \bigl/{\rm
Im}\sigma_4$, then $\mathrm{dim}\, K_2=1$. If we define a map
\begin{math}
  \tau : S^2T^*\otimes F_2 \to \BbbC
\end{math}
by 
\begin{math}
  \tau(D) = D_{11}^1 - 2D_{12}^1 - D_{22}^2 + 2D_{12}^2 ,
\end{math}
then, the sequence 
\begin{displaymath}
  0 \longrightarrow S^4T^* \otimes E_2
  \stackrel{\sigma_4}{\longrightarrow} S^2T^* \otimes F_2
  \stackrel{\tau}{\longrightarrow} K_2 \longrightarrow 0
\end{displaymath}
is exact.  We can deduce that the obstruction to the integrability of
$P_3$ is given by $\varphi_p=0$, where $\varphi: R_3 \to K_2$ is
defined by
\begin{displaymath}
  \varphi(s)_p =\; [\nabla p_0(P_3(s))]^1_{11} \; -2[\nabla
  p_0(P_3(s))]^1_{12} \; -\; 
  [\nabla p_0(P_3(s))]^2_{22}\;  +2[\nabla p_0(P_3(s))]^2_{12}  .
\end{displaymath}
Using the equations (\ref{systP2}) and (\ref{systP3}), we obtain
\begin{displaymath}
  \begin{aligned}
    \varphi(s)&=\nabla_{11}[2s_{21}-s_{22}-ss_2+2ss_1 +
    \mathcal{R}s+\mathcal{R}_2]-2\nabla_{12} [2s_{21}-s_{22} -
    ss_2+2ss_1+ \mathcal{R}s +\mathcal{R}_2]
    \\
    &+2\nabla_{12}[2s_{21}-s_{11}-2ss_2+ss_1+\mathcal{R}s +
    \mathcal{R}_1]- \nabla_{22}[2s_{21}-s_{11}-2ss_2+ss_1 +
    \mathcal{R}s + \mathcal{R}_1]
  \end{aligned}
\end{displaymath}
By the formula (\ref{permutation}) we can eliminate the 4$^{th}-$order
derivatives and find the expression of $\varphi$.$\Box$

\smallskip

We can remark that ${\rm dim} \, g_{3, p}=0$ and therefore ${\rm dim}
\, g_{k, p}=0$ for every $k >3$.  It follows that $P_3$ is involutive.

\medskip

{\sc Remark}. If ${\mathcal R} = 0$, then $\varphi = 0$, therefore, all
3$^{rd}-$order solution of $P_3$ can be lifted into a
4$^{th}-$solution.  Since $P_3$ is involutive $P_3$ is formally
integrable and consequently, it is integrable in the analytical case.
We have the following result:
    
\begin{corollaire}
  \label{R=0}
  If $\W$ is a parallelizable 3-web on the plane, then for all $L_0 \in
  E_p$ there exists a germ of linearizations $L$ which prolongs $L_0$.
\end{corollaire}

In accord of the Graf-Sauer Theorem, one can deduce that for a
parallelizable web, there are non projectively equivalent linearizations.
Indeed, it is sufficient to consider $L_0, L_0' \in E_p$ with $s_p \neq
s'_p$ and to prolong them in germs of linearization to obtain two non
projectively equivalent germs of linearization.

\bigskip

\section{Second obstruction}

In the sequel we will suppose that $\R \neq 0$.  In this case the
compatibility condition (\ref{varphi}) is not satisfied, so we have to
introduce into our differential system and consider the second order
quasi-linear system $P_{\varphi}=0$:
\begin{displaymath}
  P_{\varphi}:= (P_2,\varphi),
\end{displaymath}
where $P_2$ is defined by (\ref{systP2}) and $\varphi$ is given by the
equation (\ref{varphi}). The diagram associated to $P_{\varphi}$ is:
\begin{displaymath}
  \begin{CD}
    &&S^{3} T^*\otimes E_2 @>{\sigma _{3}(P_{\varphi})}>> (T^*\otimes
    F_2) \oplus (T^*\otimes K_2) @>{\tau_3 }>> K_3 @>>> 0
    \\
    && @VV{\varepsilon }V @VV{\varepsilon}V
    \\
    R_{3} @>>> J_{3}\, E_2 @>{p_1(P_{\varphi})}>> J_1\, F_2\oplus
    J_1K_2
    \\
    @VV{\overline{\pi}}V @VV{\pi}V @VV{\pi}V
    \\
    R_2@>>> J_2\,E_2 @>{p_0(P_{\varphi})}>> F_2\oplus K_2
  \end{CD}
\end{displaymath}

\begin{lemme}
  \label{prop} 
  A 2$^{nd}-$order formal solution $j_{2,p}s$ of $P_\varphi$ at $p\in M$, 
  can be lifted into a 3$^{rd}-$order solution if and only if:
  \begin{displaymath}
    \left\{
      \begin{aligned}
        \psi^1_p s &:=
        24\mathcal{R}s_2^2-48\mathcal{R}s_1s_2+\alpha(s)s_1 +
        \beta(s)s_2+\gamma(s) =0
        \\
        \psi^2_p s &:=
        -24\mathcal{R}s_1^2+48\mathcal{R}s_1s_2+\hat{\alpha}(s)s_1 +
        \hat{\beta}(s)s_2+\hat{\gamma}(s) =0.
      \end{aligned}
    \right.
  \end{displaymath}
  where $\alpha$, $\beta$, $\hat{\alpha}$, $\hat{\beta}$ are
  polynomials in $s$ of degree $2$ with coefficients $\mathcal{R}$ and
  its derivatives up to order $2$, $\gamma$ and $\hat{\gamma}$ are
  polynomials in $s$ of degree $3$ with coefficients ${\mathcal R}$
  and its derivatives up to order 4. Their explicit expressions are
  given in Appendix.
\end{lemme}

{\it Proof.}  The symbol of differential operator $\varphi$ is
$\sigma_2({\varphi)}:S^2T^* \otimes E_2 {\longrightarrow} K_2$ and its
prolongation $\sigma_3(\varphi):S^3T^* \otimes E_2 {\longrightarrow}
T^* \otimes K_2 $ are given by:
\begin{displaymath}
  \sigma_2(\varphi)(A) = -24{\mathcal R}A_{21} \quad \mathrm{and} \quad  
  \sigma_3(\varphi)(B)(e_i) =-24{\mathcal R}B_{i21}, \quad i=1,2
\end{displaymath}
where $A_{21}:=A(e_2,e_1)$ and $B_{i21} = B(e_i,e_2,e_1)$ are the
components of the corresponding tensors with respect to the adapted
basis $\{e_1,e_2\}$. 

\smallskip

Note, that we have
\begin{math}
  g_2(P_\varphi)=g_2(P_2)\cap {\rm Ker} \sigma_2 (\varphi) =0,
\end{math}
therefore for every $\ell>2$ we obtain that 
\begin{math}
  g_\ell(P_\varphi)=0,
\end{math}
and so $P_{\varphi}$ is involutive.

\smallskip

The kernel of the symbol of the first prolongation is
$g_3(P_\varphi)$ defined by the system
\begin{equation}
  \label{g3(Pphi)}
  \left\{
    \begin{aligned}
      A^1_1 &:= B_{122}-2B_{112}=0,
      \\
      A^1_2 &:= B_{222}-2B_{122}=0,
      \\
      A^2_1 &:= B_{111}-2B_{112}=0,
      \\
      A^2_2 &:= B_{112}-2B_{122}=0,
      \\
      C_1 &:= -24\mathcal{R} B_{112}=0,
      \\
      C_2 &:= -24\mathcal{R} B_{122}=0.
    \end{aligned}
  \right.
\end{equation}
There are two relations in this system (\ref{g3(Pphi)}). Namely
$24\mathcal{R}A_1 - 2C_1 + C_2 = 0$ and $24\mathcal{R}B_2 + C_1 - 2C_2
= 0$.  So ${\rm rank}\,\sigma_3(P_\varphi)=4$, and
\begin{displaymath}
  {\rm dim} \, K_3 = {\rm dim } \bigl( (T^*\otimes F_2)\oplus
  (T^*\otimes K_2) \bigl/ {\rm Im}\sigma_3\bigl) = 2.
\end{displaymath}
Moreover, if we define
\begin{displaymath}
  \begin{array}{lclc} 
    \tau_3 \, : \, &(T^*\otimes F_2)\oplus
    (T^* \otimes K_2) &
    {\longrightarrow} & K_3 \simeq \BbbC^2\\
    &{ (A^1,\, A^2, \, C)}& \longmapsto& {(D^1,D^2)}
  \end{array}
\end{displaymath}
by:
\begin{displaymath}
  D^1 := 24\mathcal{R}A_1 - 2C_1 + C_2, \qquad \mathrm{and} \qquad
  D^2 := 24\mathcal{R}B_2 + C_1 - 2C_2,
\end{displaymath}
then the sequence
\begin{displaymath}
  0 \longrightarrow  S^3T^*\otimes E_2 \stackrel{\sigma_\varphi }{\longrightarrow}  
  (T^*\otimes F_2)\oplus (T^*\otimes K_2)  
  \stackrel{\tau_3} {\longrightarrow}  K_3 \longrightarrow 0  
\end{displaymath}
is exact.  We can deduce that a 2$^{nd}$ order solution $(j_2s)_p$ of
$P_\varphi$ can be lifted into a 3$^{rd}$ order solution if and only if
$[\tau_3\nabla (P_\varphi(s))]_p=0$. Let
\begin{displaymath}
  (\psi^1, \psi^2)_p:=[\tau_3\nabla (P_\varphi(s))]_p=\tau_3(\nabla P_2(s) , \,
  \nabla\varphi)_p
\end{displaymath}
We have:
\begin{displaymath}
  \begin{aligned}
    \psi^1 &= 24\mathcal{R}[\N(P_2(s))]_1 +[\nabla(\varphi)]_2
    -2[\nabla(\varphi)]_1
    \\
    \psi_2 & = 24\mathcal{R}[\N(P_2(s))]_2 +[\nabla(\varphi)]_1
    -2[\nabla(\varphi)]_2
  \end{aligned}
\end{displaymath}
Using the equations $P_2(s)_p=0$ and $\varphi (s)_p =0$ and the permutation
formula (\ref{permutation}), we find that $\psi^1$ and $\psi^2$ can be
written as a function of $s$ and its derivatives up to order 3.
Nevertheless, using the formula (\ref{permutation}) we can also eliminate
the 3$^{th}$ order derivatives of $s$. On the other hand, with the help of
the equation $P_2=0$ and $\varphi=0$ we can express the 2$^{nd}$ order
derivatives of $s$ too with the 1$^{st}$ order derivatives of $s$.  The
calculation carried out with MAPLE gives the formulas.

\bigskip

\section{The linearization theorem}

Since the compatibility conditions $\psi^1=0$ and $\psi^2=0$ found in
the previous section are not identically satisfied, we have to
introduce them into the system $P_{\varphi}$. We arrive at the system:
\begin{displaymath}
  P_\psi = (P_2, \varphi, \psi^1, \psi^2).
\end{displaymath}
Differentiating the equations $\psi^1 = 0$ and $\psi^2= 0$ with
respect to $e_1$ and $e_2$ we find 4 equations:
\begin{displaymath}
  {\small
    \left\{
      \begin{aligned}
        \psi^1_1 = &
        24\mathcal{R}_1s_2^2+(48\mathcal{R}s_2-48\mathcal{R}s_1 +
        \beta)s_{12} -48\mathcal{R}_1s_1s_2+ (\alpha
        -48\mathcal{R}s_2)s_{11} + \alpha_1s_1 + \beta_1s_2+\gamma_1
        \\
        \psi^1_2 = & 24\mathcal{R}_2s_2^2+(48\mathcal{R}s_2 -
        48\mathcal{R}s_1+ \beta)s_{22} -48\mathcal{R}_2s_1s_2 +
        (\alpha-48\mathcal{R}s_2)s_{21}+\alpha_2s_1+\beta_2s_2+\gamma_2
        \\
        \psi^2_1 = & -24\mathcal{R}_1s_1^2+ \hat{\alpha}s_{11} -
        48\mathcal{R}s_1+48\mathcal{R}s_2 +48\mathcal{R}_1s_1s_2 +
        (48\mathcal{R}s_1+\hat{\beta})s_{12} + \hat{\alpha}_1s_1 +
        \hat{\beta}_1s_2+\hat{\gamma}_1
        \\
        \psi^2_2 = & -24\mathcal{R}_2s_1^2 + ( 48\mathcal{R}s_2
        -48\mathcal{R}s_1 + \hat{\alpha})s_{21}
        +48\mathcal{R}_2s_1s_2+ (48\mathcal{R}s_1 +\hat{\beta})s_{22}
        +\hat{\alpha}_2s_1 +\hat{\beta}_2s_2+\hat{\gamma}_2
      \end{aligned}
    \right.}
\end{displaymath}
In this expression, we can eliminate the second order derivatives using the
equation $P_2 = 0$ and $\varphi = 0$, and with the help of the equation
$\psi^1 = 0$ and $\psi^2 = 0$, we can express the terms $s_1^2$ and $s_2^2$
as a function of $s_1$, $s_2$ and the product $s_1s_2$.  Therefore the
system
\begin{displaymath}
  P_\psi = 0, \ \ \nabla P_{\psi^1} = 0,\ \  \nabla P_{\psi^2} = 0
\end{displaymath}
is equivalent to the system formed by the equation $P_\psi=0$ and the
four {\it linear} equations in $s_1$, $s_2$ and $s_1s_2$:
\begin{eqnarray}
  \label{abc}
  {\mathcal S} = \left\{
    \begin{aligned}{}
      a^1s_1+b^1s_2+c^1s_1s_2=d^1,
      \\
      a^2s_1+b^2s_2+c^2s_1s_2=d^2,
      \\
      a^3s_1+b^3s_2+c^3s_1s_2=d^3,
      \\
      a^4s_1+b^4s_2+c^4s_1s_2=d^4,
    \end{aligned}
  \right.
\end{eqnarray}
where $a^i$, $b^i$, $i=1, ...,4$ are polynomials in $s$ of degree $3$,
whose coefficients are $\mathcal{R}$ and its derivatives up to order
3, $c^1$ and $c^4$ are polynomials in $s$ of degree $1$ with
coefficients $\R$, $\R_1$ and $\R_2$, $c^2$ and $c^3$ can be expressed
as a function of $\R$, $\R_1$ and $\R_2$, and $d^1$, $d^4$ (resp.
$d^2$ and $d^3$) are polynomials in $s$ of degree $5$ (resp.  $4$),
with coefficients $\R$ and its derivatives up to order 5.  Its
explicit expressions will be given in Appendix.

\bigskip

The direct computation by MAPLE shows us  that the determinant
\begin{displaymath}
  \left|
    \begin{array}{llll}
      a^1&b^1&c^1&d^1\\
      a^2&b^2&c^2&d^2\\
      a^3&b^3&c^3&d^3\\
      a^4&b^4&c^4&d^4
    \end{array}
  \right| 
\end{displaymath}
is identically null, so that the system ${\mathcal S}$ is compatible.
On the other hand, the 3$^{rd}-$order minors of the system ${\mathcal
  S}$ are polynomials in $s$ of degree 7 which are not identically
zero. There is a open dense ${\mathcal U} \subset \BbbC^2$ on which,
\begin{displaymath}
  D(s) = 
  \left|
    \begin{array}{lll}
      a^1&b^1 &c^1\\
      a^2&b^2&c^2 \\
      a^3&b^3&c^3
    \end{array}
  \right| \neq 0.
\end{displaymath}
Solving on ${\mathcal U}$ the system ${\mathcal S}$ for $s_1$, $s_2$
and $s_1s_2$ we obtain:
\begin{equation}
  \label{s1s2}
  s_1 = F(s) = \dd\frac{A(s)}{D(s)},
  \qquad
  s_2= G(s) = \dd\frac{B(s)}{D(s)}
\end{equation}
and
\begin{equation}
  \label{produit}
  s_1s_2 = H(s) = \dd\frac{C(s)}{D(s)},
\end{equation}
where the polynomials in  $s$:
\begin{displaymath}
  A(s) = 
  \left|
    \begin{array}{rrr}
      -d^1&b^1&c^1\\
      -d^2&b^2&c^2\\
      -d^3&b^3&c^3
    \end{array}
  \right
  | , \quad 
  B(s) =
  \left|
    \begin{array}{rrr}
      a^1&-d^1&c^1\\
      a^2&-d^2&c^2\\
      a^3&-d^3&c^3
    \end{array}
  \right| , \quad
  C(s) =
  \left|
    \begin{array}{rrr}
      a^1&b^1&-d^1\\
      a^2&b^2&-d^2\\
      a^3&b^3&-d^3
    \end{array}
  \right|
\end{displaymath}
are of degrees $8$, $8$, and $11$ respectively.
    
By (\ref{produit}) we must find $F(s)\, G(s) = H (s)$. Thus, the
solution of $s$ for the linearization system must be in the algebraic
manifold defined by
\begin{equation}
  Q_1(s) : =  AB -CD = 0.
\end{equation} 
On the other hand, the compatibility condition of the system
(\ref{s1s2}) is $s_{12} - s_{21} = \R s$.  Computing it explicitly we
find that $s$ must be in the algebraic manifold defined by
\begin{displaymath}
  Q_2 (s)=0,
\end{displaymath}
where $Q_2$ is polynomial in $s$ of degree 15. Indeed, if
\begin{math}
  A(s) = \sum_{i=1}^8 A_i s^i, 
\end{math}
\begin{math}
  B(s) = \sum_{i=1}^8 B_i s^i,
\end{math}
and
\begin{math}
  D(s) = \sum_{i=1}^7 D_i s^i
\end{math}
where $A_i$, $B_i$ and $C_i$ are function on $M$, then using
(\ref{s1s2}) we obtain
\begin{alignat*}{1}
  Q_2(s) & = \Bigl( \sum_{i=1}^8 (\nabla_2 B_i) s^i \Bigl)
  \Bigl( \sum_{i=1}^7 D_i s^i \Bigl) 
  - \Bigl( \sum_{i=1}^8 B_i s^i \Bigl) \Bigl( \sum_{i=1}^7 (\nabla_2
    D_i) s^i \Bigl) 
  \\& 
  - \Bigl( \sum_{i=1}^8 (\nabla_1 A_i) s^i \Bigl)
  \Bigl( \sum_{i=1}^7 D_i s^i \Bigl) 
  - \Bigl( \sum_{i=1}^8 A_i s^i \Bigl) \Bigl( \sum_{i=1}^7 (\nabla_1
    D_i) s^i \Bigl)
  \\
  & + \Bigl( \sum_{i=1}^8 B_i s^{i-1} \Bigl)\Bigl( \sum_{i=1}^8 A_i
    s^i \Bigl)
    - \Bigl( \sum_{i=1}^8 B_i s^{i} \Bigl)\Bigl( \sum_{i=1}^8 A_i
    s^{i-1} \Bigl) - \R s D^2.
\end{alignat*}
Moreover, we must impose that $s_1$ and $s_2$ given by (\ref{s1s2})
verify the 5 equations of $P_\psi$, this implies $5$ polynomial
equations $Q_i = 0$, $i=3,...,7$.  Finally, we arrive at the
conclusion that if the web is linearizable then $s$ must be in the
algebraic manifold ${\mathcal A}$, where ${\mathcal A}$ is defined by
the equations $Q_i = 0$, $i=1,...,7$:
\begin{displaymath}
  \mathcal{A}:=\{ Q_i=0 \ | \ i=1,...,7\}.
\end{displaymath}

\smallskip

So the compatibility system (therefor the linearization system) has a
solution in the neighborhood of a point $p \in M$ if and only if the
algebraic variety ${\mathcal A}$ is not empty.  If ${\mathcal A}\neq
\emptyset $, then for all smooth point $s_0 \in {\mathcal A}$, there
exists a neighborhood $U$ of $s_0$ so that all $s\in U$ can be
prolonged in a germ $\tilde{s}$ as a basis of linearization.  The
explicit expression of the polynomials $Q_i$ can be computed with the
help of MAPLE.  The degree of these polynomials $Q_i$, $i=1 ...  7$
are 18, 15, 23, 23, 24, 17 and 17 respectively.  One obtains the
following results:

\medskip
                               
\begin{theoreme}  
  A non-parallelizable $3-$web $\W$ is linearizable if and only if
  there is an open set $U$ of $M$ on which the polynomials $Q_1$, ...,
  $Q_7$ have common zeros.  Moreover, if this condition is satisfied,
  then for all $p\in U$ and all pre-linearization $L_0\in E_p$ whose
  base is in ${\mathcal A}=\{ Q_i=0 \, | \, i=1,...,7\}$, there exists
  a unique linearization $L$ so that $L_p = L_0$.
\end{theoreme}
Since the lowest degree of the polynomials defining ${\mathcal A}$ is 15
we arrive at the
\begin{theoreme}
  For a non parallelizable $3-$web, there exists at most $15$
  projectively non equivalent linearizations.
\end{theoreme}

\noindent
Finally, the \textbf{Gronwall conjecture} can be expressed now in the
following way: for any non parallelizable 3-web on a 2-dimensional
manifold
\begin{displaymath}
  deg[Rad(Q_1,...,Q_7)]=1,
\end{displaymath}
where $Rad$ denotes the radical of the corresponding polynomials and
$deg$ is its degree.

\bigskip

\noindent
{\bf Examples}

\begin{enumerate}
  {\small
  \item Consider the web $\W$ defined by $x=cte$, $y=cte$,
    $f(x,y)=cte$, where $f(x,y):=(x+y)e^{-x}$.  This web is not
    parallelizable in a neighborhood of $(0,0)$ because the Chern
    connection is not flat.  Indeed, the component of the curvature
    tensor ${\mathcal R}$ at $(0,0)$ can be computed directly from the
    function $f$ by the formula
    \begin{displaymath}
      {\mathcal R}=\frac{1}{f_x f_y} \Bigl(\frac{f_{xxy}}{f_x} -
      \frac{f_{xyy}}{f_y} + \frac{f_{xy} f_{yy}}{f_y^2} - \frac{f_{xx}
        f_{xy}}{f_x^2} \Bigl),
    \end{displaymath}
    (cf. \cite{AS}, p. 24). In this example we have $\mathcal R_{(0,0)}=-1$.  The
    computation gives that $Rad(Q_1,\cdots, Q_7)=s+1$ on a neighborhood of
    $(0,0)$. Thus the web { is linearizable} in a neighborhood of $(0,0)$
    and all the linearizations are projectively equivalent.

\smallskip

\item Let $\W$ be the web defined by $x=cte$, $y=cte$, $f(x,y)=cte$,
  where $$f= \log(x)+ \tfrac{1}{2}\log \left( \frac{x^2+y^2}{x^2}
  \right) + {\rm arctg} \left( \frac{y}{x} \right).$$
  We have
  ${\mathcal R}_{(1,0)}=2$, so $\mathcal{W}$ is not parallelizable at
  $(1,0).$ On the other hand the resultant of the polynomials $Q_2,
  Q_6$ is not zero at $(1,0)$. So this web { is not linearizable} at
  $(1,0)$.}

\end{enumerate}

\section{Appendix}

\footnotesize

\noindent
\begin{alignat*}{1}
  \alpha & = 30 { R} s^2-18 { R}2 s-\tfrac{3}{4{ R}} (-16 { R} {
    R}_{22}+14 { R}_2^2-40 { R}^3-56 { R}_1 { R}_2+40 { R} { R}_{12}
  +56 { R}_1^2-40 { R} { R}_{11}),
  \\
  \beta & =\tfrac{3}{4{ R}} (24 { R} { R}_2-24 { R} { R}_1) s -15 { R}
  s^2 -\tfrac{3}{4{ R}} (70 { R}_1 { R}_2-44 { R} { R}_{12}+20 { R} {
    R}_{11}+20 { R} { R}_{22}-28 { R}_1^2-28 { R}_2^2-60 { R}^3),
  \\
  \gamma & =-\tfrac{3}{4{ R}} (-6 { R} { R}_1+3 { R} { R}_2)
  s^3-\tfrac{3}{4{ R}} (7 { R}_2 { R}_{12}+12 { R} { R}_{112}-14 {
    R}_1 { R}_{12}-7 { R}_2 { R}_{22}+14 { R}_1 { R}_{11}-8 { R} {
    R}_{111}
  \\
  & +4 { R} { R}_{222}-47 { R}_2 { R}^2+14 { R}_1 { R}_{22}-7 {
    R}_{11} { R}_2-30 { R}_1 { R}^2 -12 { R} { R}_{122})
  s-\tfrac{3}{4{ R}} (-7 { R}_2 { R}_{112}+7 { R}_2 { R}_{122}
  \\
  & +35 { R}_1 { R}_2 { R}-38 { R}_1^2 { R}-2 { R}_2^2 { R}+12 {
    R}_{1122} { R}-8 { R}_{1112} { R}-4 { R}_{1222} { R}-48 { R}^2 {
    R}_{12}
  \\
  & +8 { R}^2 { R}_{11}+40 { R}^2 { R}_{22}+8 { R}^4 -14 { R}_1 {
    R}_{122}+14 { R}_1 { R}_{112}).
  \\
  \hat{\alpha} & = \tfrac{3}{4{ R}} (24 { R} { R}_2-24 { R} { R}_1) s
  -15 { R} s^2 -\tfrac{3}{4{ R}} ( 20 { R} { R}_{22}+70 { R}_1 {
    R}_2-28 { R}_1^2+60 { R}^3-28 { R}_2^2-44 { R} { R}_{21}+20 { R} {
    R}_{11}),
  \\
  \hat{\beta} & = 30 { R} s^2+18 { R}_1 s-\tfrac{3}{4{ R}} (40 {
    R}^3-16 { R} { R}_{11}-56 { R}_1 { R}_2+56 { R}_2^2+14 { R}_1^2+40
  { R} { R}_{21}-40 { R} { R}_{22}),
  \\
  \hat{\gamma} & =-\tfrac{3}{4{ R}} (3 { R} { R}_1-6 { R} { R}_2)
  s^3-\tfrac{3}{4{ R}} (-7 { R}_1 { R}_{22}-1 4 { R}_2 { R}_21+12 { R}
  { R}_{221}+14 { R}_{11} { R}_2+7 { R}_1 { R}_{21}-12 { R} { R}_{211}
  \\
  & - 8 { R} { R}_{222}+14 { R}_2 { R}_{22}-7 { R}_1 { R}_{11}+4 { R}
  { R}_{111}+30 { R}_2 { R}^2+47 { R}_1 { R}^2) s-\tfrac{3}{4{ R}} (35
  { R}_1 { R}_2 { R}-7 { R}_1 { R}_{211}
  \\
  & +7 { R}_1 { R}_{221} +8 { R}_{2221} { R}+4 { R}_{2111} { R}-12 {
    R}_{2211} { R}+8 { R}^2 { R}_{22}-2 { R}_1^2 { R}+40 { R}^2 {
    R}_{11}-48 { R}^2 { R}_{21}-38 { R}_2^2 { R}
  \\
  &-8 { R}^4+14 { R}_2 { R}_{211}-14 { R}_2 { R}_{221}),
\end{alignat*}

\begin{alignat*}{1}
  a^1 & = \tfrac{363}{2} { R}_{12} -\tfrac{297}{2} { R}
  s^3+(\tfrac{441}{4} { R}_2-72 { R}_1) s^2 +(\tfrac{825}{4} {
    R}_2-102 { R}_1) { R}+\tfrac{1}{{ R}}((-36 { R}_{22}-147 {
    R}_{11}+114 { R}_{12}) { R}_1
  \\
  & +(\tfrac{1}{{ R}} (\tfrac{381}{2} { R}_1^2-\tfrac{1023}{4} { R}_1
  { R}_2+96 { R}_2^2) -\tfrac{1305}{2} { R}^2-\tfrac{273}{2} {
    R}_{11}-\tfrac{165}{2} { R}_{22} ) s + 15 { R}_{122}-33 {
    R}_{112}-3 { R}_{222}
  \\
  & +36 { R}_{111} +(\tfrac{57}{4} { R}_{22}+ \tfrac{189}{4} {
    R}_{11}-\tfrac{177}{4} { R}_{12}) { R}_2)+\tfrac{1}{{
      R}^2}(\tfrac{231}{2} { R}_1^3-\tfrac{525}{4} { R}_1^2 {
    R}_2+\tfrac{273}{4} { R}_1 { R}_2^2-\tfrac{63}{4} { R}_2^3),
  \\
  a^2 & = -\tfrac{9}{2} { R} s^3+(-\tfrac{45}{2} { R}_1+18 { R}_2)
  s^2+(\tfrac{39}{2} { R}^2-3 { R}_{22}-\tfrac{39}{2} {
    R}_{12}+\tfrac{39}{2} { R}_{11}+ \tfrac{1}{{ R}}(-\tfrac{57}{2} {
    R}_1^2+42 { R}_1 { R}_2-\tfrac{21}{8} { R}_2^2)) s
  \\
  &(-\tfrac{342}{2} { R}_1+ \tfrac{417}{2} { R}_2) { R}-42 {
    R}_{122}+12 { R}_{222}+42 { R}_{112}+\tfrac{1}{{ R}}((24 {
    R}_{22}-\tfrac{45}{2} { R}_{11}-\tfrac{87}{2} { R}_{12}) { R}_1+
  (-\tfrac{39}{2} { R}_{22}
  \\
  & -\tfrac{15}{2}{ R}_{11}+\tfrac{81}{2} { R}_{12}) {
    R}_2)+\tfrac{1}{{ R}^2}(\tfrac{63}{2} { R}_1^3-\tfrac{21}{2} {
    R}_1^2 { R}_2-\tfrac{105}{8} { R}_1 { R}_2^2+\tfrac{21}{4} {
    R}_2^3),
  \\
  a^3 & = -\tfrac{9}{4} { R} s^3+(-18 { R}_1+\tfrac{117}{4}{ R}_2)
  s^2+(\tfrac{999}{4} { R}^2-\tfrac{129}{4} { R}_{12}+\tfrac{39}{4} {
    R}_{22}+\tfrac{39}{4} { R}_{11}+\tfrac{1}{{ R}}(-\tfrac{57}{4} {
    R}_1^2+\tfrac{429}{8}{ R}_1 { R}_2
  \\
  & -\tfrac{111}{4} { R}_2^2)) s+(-\tfrac{429}{2} {
    R}_1+\tfrac{159}{4} { R}_2) { R}+ 48 { R}_{112}+6 { R}_{222}-30 {
    R}_{122}-18 { R}_{111}+\tfrac{1}{{ R}}((\tfrac{9}{2} {
    R}_{22}+\tfrac{75}{2} { R}_{11}
  \\
  & -\tfrac{141}{2} { R}_{12}) { R}_1+(\tfrac{39}{4} {
    R}_{22}-\tfrac{45}{4} { R}_{11}+\tfrac{39}{4} { R}_{12}) {
    R}_2)+\tfrac{1}{{ R}^2}(-\tfrac{21}{2} { R}_1^3+\tfrac{21}{2} {
    R}_1^2 { R}_2+\tfrac{231}{8} { R}_1 { R}_2^2-\tfrac{63}{4} {
    R}_2^3),
  \end{alignat*}

  \begin{alignat*}{1}
  b^1&= q144 { R} s^3+(-63 { R}_2+\tfrac{531}{4} { R}_1) s^2+(120 {
    R}^2+156 { R}_{22}-156 { R}_{12}+57 { R}_{11}+\tfrac{1}{{
      R}}(-\tfrac{183}{4} { R}_1^2+\tfrac{465}{2} { R}_1 { R}_2
  \\
  & -219 { R}_2^2)) s+(-99 { R}_2+\tfrac{279}{4} { R}_1) { R}+39 {
    R}_{112}-21 { R}_{122}-15 { R}_{111}+\tfrac{1}{{
      R}}((\tfrac{291}{4} { R}_{11} +\tfrac{159}{4} {
    R}_{22}-\tfrac{423}{4} { R}_{12}) { R}_1
  \\
  & +(60 { R}_{12}-75 { R}_{11}-27 { R}_{22}) { R}_2)+\tfrac{1}{{
      R}^2}(-\tfrac{231}{4} { R}_1^3+\tfrac{609}{4} { R}_1^2 {
    R}_2-\tfrac{357}{4} { R}_1 { R}_2^2+\tfrac{63}{2} { R}_2^3),
  \\
  b^2&= \tfrac{9}{4} { R} s^3+(\tfrac{117}{4}{ R}_1-18 { R}_2)
  s^2+(\tfrac{741}{4} { R}^2-3\tfrac{9}{4} { R}_{11}-3\tfrac{9}{4} {
    R}_{22}+12\tfrac{9}{4} { R}_{12}+\tfrac{1}{{ R}}(\tfrac{111}{4} {
    R}_1^2-\tfrac{429}{8}{ R}_1 { R}_2+\tfrac{57}{4} { R}_2^2)) s
  \\
  & +(\tfrac{603}{4} { R}_1+\tfrac{39}{2} { R}_2) { R}+48 {
    R}_{122}-30 { R}_{112}-18 { R}_{222}+6 { R}_{111}+ \tfrac{1}{{
      R}}((\tfrac{39}{4} { R}_{11}+3\tfrac{9}{4} {
    R}_{12}-\tfrac{45}{4} { R}_{22}) { R}_1
  \\
  & +(\tfrac{9}{2} { R}_{11}-\tfrac{141}{2} { R}_{12}+\tfrac{75}{2} {
    R}_{22}) { R}_2)+\tfrac{1}{{ R}^2}(-\tfrac{63}{4} {
    R}_1^3+\tfrac{231}{8} { R}_1^2 { R}_2+\tfrac{21}{2} { R}_1 {
    R}_2^2-\tfrac{21}{2} { R}_2^3),
  \\
  b^3&= \tfrac{9}{2} { R} s^3+(-\tfrac{45}{2} { R}_2+18 { R}_1)
  s^2+(-\tfrac{39}{2} { R}^2+3 { R}_{11}+\tfrac{39}{2} {
    R}_{12}-\tfrac{39}{2} { R}_{22}+ \tfrac{1}{{ R}}(\tfrac{21}{8} {
    R}_1^2-42 { R}_1 { R}_2+\tfrac{57}{2} { R}_2^2)) s
  \\
  & +(\tfrac{243}{2}{ R}_2+\tfrac{9}{2} { R}_1) { R}-42 { R}_{112}+42
  { R}_{122}+12 { R}_{111}+\tfrac{1}{{ R}}((-\tfrac{39}{2} {
    R}_{11}+\tfrac{81}{2} { R}_{12}-\tfrac{15}{2}{ R}_{22}) { R}_1+(24
  { R}_{11}
  \\
  & -\tfrac{87}{2} { R}_{12}-\tfrac{45}{2} { R}_{22}) {
    R}_2)+\tfrac{1}{{ R}^2}(\tfrac{21}{4} { R}_1^3-\tfrac{105}{8} {
    R}_1^2 { R}_2-\tfrac{21}{2} { R}_1 { R}_2^2+\tfrac{63}{2} {
    R}_2^3),
    \\
    \\
    c^1&= (234 { R} s+18 { R}_2+18 { R}_1)
    \\
    c^2&= (36 { R}_1-18 { R}_2),
    \\
    c^3&= (18 { R}_1-36 { R}_2)
  \end{alignat*}
 
  \begin{alignat*}{1}
    d^1 & = \tfrac{45}{8} { R} s^5+(-\tfrac{171}{8} {
      R}_1+\tfrac{45}{2} { R}_2) s^4+(\tfrac{9}{2} {
      R}_{11}-\tfrac{9}{2} { R}_{22}+\tfrac{1}{{ R}}(-\tfrac{99}{8} {
      R}_1^2+\tfrac{63}{8} { R}_1 { R}_2-\tfrac{27}{8} { R}_2^2))
    s^3+((-\tfrac{21}{4} { R}_2
    \\
    &-411 { R}_1) { R}-\tfrac{423}{8} { R}_{122}+\tfrac{423}{8} {
      R}_{112}+33 { R}_{222}-\tfrac{51}{2} { R}_{111}+\tfrac{1}{{
        R}}((-\tfrac{375}{8} { R}_{12}+\tfrac{375}{8} {
      R}_{22}+\tfrac{375}{8} { R}_{11}) { R}_1
    \\
    & +(-\tfrac{111}{2} { R}_{11}-\tfrac{111}{2} {
      R}_{22}+\tfrac{111}{2} { R}_{12}) { R}_2)) s^2+(-\tfrac{2205}{8}
    { R}^3+(\tfrac{1233}{4} { R}_{22}-\tfrac{501}{4} {
      R}_{12}+\tfrac{87}{4} { R}_{11} -\tfrac{567}{2} { R}_{21}) { R}
    \\
    & -\tfrac{363}{2} { R}_1^2+\tfrac{903}{8} { R}_1 {
      R}_2-\tfrac{417}{8} { R}_2^2-\tfrac{69}{2} { R}_{1112}-36 {
      R}_{1222}+\tfrac{135}{2} { R}_{1122}+6 { R}_{1111} +\tfrac{1}{{
        R}}((-\tfrac{567}{8} { R}_{122}
    \\
    & +\tfrac{567}{8} { R}_{112}+\tfrac{15}{2} {
      R}_{222}-\tfrac{33}{2} { R}_{111}) { R}_1+(-\tfrac{9}{2} {
      R}_{222}+63 { R}_{122} -63 { R}_{112}+3 { R}_{111}) {
      R}_2-\tfrac{129}{8} { R}_{11}^2
    \\
    & +(\tfrac{99}{4} { R}_{12}-\tfrac{69}{4} { R}_{22}) { R}_{11}
    +\tfrac{39}{4} { R}_{12} { R}_{22}-\tfrac{9}{8} { R}_{22}^2
    -\tfrac{69}{8} { R}_{12}^2)+\tfrac{1}{{ R}^2}((\tfrac{231}{8} {
      R}_{22}-\tfrac{231}{8} { R}_{12}+\tfrac{231}{8} { R}_{11}) {
      R}_1^2
    \\
    & +(-\tfrac{147}{8} { R}_{22}-\tfrac{147}{8} {
      R}_{11}+\tfrac{147}{8} { R}_{12}) { R}_2 { R}_1+(\tfrac{63}{8} {
      R}_{22}-\tfrac{63}{8} { R}_{12}+\tfrac{63}{8} { R}_{11}) {
      R}_2^2)) s+(\tfrac{39}{4} { R}_2-\tfrac{303}{8} { R}_1) { R}^2
    \\
    & +(-6 { R}_{111}+\tfrac{138}{8} { R}_{112}-\tfrac{135}{8} {
      R}_{122}) { R}+(\tfrac{465}{8} { R}_{22} +\tfrac{741}{8} {
      R}_{11}-\tfrac{549}{8} { R}_{12}-63 { R}_{21}) {
      R}_1+(-\tfrac{207}{4} { R}_{22}
    \\
    & + 21 { R}_{12}-\tfrac{231 }{4} { R}_{11}+\tfrac{315}{4} {
      R}_{21}) { R}_2+6 { R}_{11112}-9 { R}_{11122}+3 {
      R}_{11222}+\tfrac{1}{{ R}}(-\tfrac{627}{8} {
      R}_1^3+\tfrac{717}{8} { R}_1^2 { R}_2
    \\
    & +(-\tfrac{33}{2} { R}_{1112}+24 { R}_{1122}-\tfrac{33}{8} {
      R}_2^2 -\tfrac{15}{2} { R}_{1222}) { R}_1+\tfrac{9}{4} {
      R}_2^3+(3 { R}_{1112}+ \tfrac{9}{2} { R}_{1222}-\tfrac{15}{2} {
      R}_{1122}) { R}_2
    \\
    & +(-\tfrac{129}{8} { R}_{112}+\tfrac{129}{8} { R}_{122}) {
      R}_{11} +(-\tfrac{69}{8}{ R}_{122}+\tfrac{69}{8} { R}_{112}) {
      R}_{12}+(-\tfrac{9}{8} { R}_{112}+\tfrac{9}{8} { R}_{122}) {
      R}_{22})
    \\
    & +\tfrac{1}{{ R}^2}((-\tfrac{231}{8} { R}_{122}+\tfrac{231}{8} {
      R}_{112}) { R}_1^2+ (\tfrac{147}{8} { R}_{122}-\tfrac{147}{8} {
      R}_{112}) { R}_2 { R}_1+(\tfrac{63}{8} { R}_{112}-\tfrac{63}{8}
    { R}_{122}) { R}_2^2),
    \end{alignat*}

    \begin{alignat*}{1}
    d^2 & = (\tfrac{9}{8}{R}_1 - \tfrac{9}{16} { R}_2) s^4+(-9 {
      R}^2-\tfrac{9}{2} { R}_{22}+9 { R}_{12}+ \tfrac{1}{{
        R}}(-\tfrac{27}{8} { R}_1^2-\tfrac{9}{16} { R}_1 {
      R}_2+\tfrac{9}{8} { R}_2^2)) s^3+((\tfrac{261}{8} {
      R}_1+\tfrac{573}{16} { R}_2) { R}
    \\
    & -\tfrac{45}{4} { R}_{112}-15/4 { R}_{222}+\tfrac{45}{4} {
      R}_{122}+ \tfrac{15}{2}{ R}_{111}+\tfrac{1}{{
        R}}((-\tfrac{69}{8} { R}_{11}+\tfrac{69}{8} {
      R}_{12}-\tfrac{69}{8} { R}_{22}) { R}_1 +(\tfrac{69}{16} {
      R}_{11}+\tfrac{69}{16} { R}_{22}
    \\
    & -\tfrac{69}{16} { R}_{12}) { R}_2)) s^2+(-\tfrac{165}{2} {
      R}^3+(135 { R}_{12}-\tfrac{165}{2} { R}_{11}-36 { R}_{22}) {
      R}-3\tfrac{45}{4} { R}_1^2+\tfrac{261}{4} { R}_1 {
      R}_2+\tfrac{51}{4} { R}_2^2 -3 { R}_{2222}
    \\
    & +\tfrac{27}{2} { R}_{1112}+\tfrac{51}{4} {
      R}_{1222}-\tfrac{81}{4}{ R}_{1122}+\tfrac{1}{{ R}}((\tfrac{3}{4}
    { R}_{222}+\tfrac{27}{8} { R}_{122}-\tfrac{27}{8} {
      R}_{112}-\tfrac{9}{2} { R}_{111}) { R}_1+(\tfrac{3}{2}{
      R}_{222}-
    \\
    & \tfrac{117}{16}{ R}_{122}+\tfrac{117}{16}{ R}_{112}-
    \tfrac{3}{2}{ R}_{111}) { R}_2+(\tfrac{15}{2}{ R}_{22}-15 {
      R}_{12}) { R}_{11}+15 { R}_{12}^2-\tfrac{45}{2} { R}_{12} {
      R}_{22}+\tfrac{15}{2}{ R}_{22}^2)
    \\
    & +\tfrac{1}{{ R}^2}((\tfrac{63}{8} { R}_{22}-\tfrac{63}{8} {
      R}_{12}+\tfrac{63}{8} { R}_{11}) { R}_1^2+ (\tfrac{21}{16} {
      R}_{22}-\tfrac{21}{16} { R}_{12}+\tfrac{21}{16} { R}_{11}) {
      R}_2 { R}_1+(\tfrac{21}{8} { R}_{12}-\tfrac{21}{8} { R}_{11}
    \\
    & -\tfrac{21}{8} { R}_{22}) { R}_2^2)) s+(12\tfrac{3}{2}{ R}_1-174
    { R}_2) { R}^2+ (111 { R}_{122}-72 { R}_{112}-45 { R}_{222}) {
      R}+(-96 { R}_{12}+78 { R}_{22}
    \\
    & +\tfrac{9}{2} { R}_{11}) { R}_1+(-\tfrac{75}{2} { R}_{11}+42 {
      R}_{12}-45 { R}_{22}) { R}_2-9 { R}_{11222}+6 { R}_{11122}+3 {
      R}_{12222}+\tfrac{1}{{ R}}(-\tfrac{171}{8} { R}_1^3
    +\tfrac{303}{16}{ R}_1^2 { R}_2
    \\
    & +(-\tfrac{3}{4} { R}_{1222}-\tfrac{9}{2} {
      R}_{1112}+\tfrac{21}{4} { R}_{1122}+3 { R}_2^2) { R}_1-
    \tfrac{3}{4} { R}_2^3+(3 { R}_{1122}-\tfrac{3}{2}{
      R}_{1222}-\tfrac{3}{2}{ R}_{1112}) { R}_2+(-15 { R}_{112}
    \\
    & +15 { R}_{122}) { R}_{12}+(-\tfrac{15}{2}{
      R}_{122}+\tfrac{15}{2}{ R}_{112}) { R}_{22})+\tfrac{1}{{ R}^2}
    ((-\tfrac{63}{8} { R}_{122}+\tfrac{63}{8} { R}_{112}) {
      R}_1^2+(\tfrac{21}{16} { R}_{112}
    \\
    & -\tfrac{21}{16} { R}_{122}) { R}_2 { R}_1+(\tfrac{21}{8} {
      R}_{122}- \tfrac{21}{8} { R}_{112}) { R}_2^2),
    \end{alignat*}

    \begin{alignat*}{1}
    d^3 & = (-\tfrac{9}{8} { R}_2+\tfrac{9}{16} { R}_1) s^4+(-9 {
      R}^2-\tfrac{9}{2} { R}_{11}+9 { R}_{12}+\tfrac{1}{{
        R}}(\tfrac{9}{8} { R}_1^2-\tfrac{9}{16} { R}_1 {
      R}_2-\tfrac{27}{8} { R}_2^2)) s^3+ ((-\tfrac{231}{8} {
      R}_2+\tfrac{1695}{16} { R}_1) { R}
    \\
    & -\tfrac{15}{2}{ R}_{222}+\tfrac{45}{4} { R}_{122}-\tfrac{45}{4}
    { R}_{112} +\tfrac{15}{4}{ R}_{111}+\tfrac{1}{{
        R}}((-\tfrac{69}{16} { R}_{22}+6\tfrac{9}{16} {
      R}_{12}-6\tfrac{9}{16} { R}_{11}) { R}_1+(-6\tfrac{9}{8} {
      R}_{12}+6\tfrac{9}{8} { R}_{22}
    \\
    & +6 \tfrac{9}{8} { R}_{11}) { R}_2)) s^2+(-36 {
      R}^3+(-\tfrac{237}{4} { R}_{11}+\tfrac{177}{2} { R}_{12}-36 {
      R}_{22}) { R}-\tfrac{897}{16} { R}_1 ^2+\tfrac{1041}{16} { R}_1
    { R}_2+\tfrac{249}{8} { R}_2^2-3 { R}_{1111}
    \\
    & -\tfrac{81}{4}{ R}_{1122}+\tfrac{27}{2} {
      R}_{1222}+\tfrac{51}{4} { R}_{1112}+\tfrac{1}{{ R}}(
    (-\tfrac{117}{16}{ R}_{112}-\tfrac{3}{2}{
      R}_{222}+\tfrac{117}{16}{ R}_{122}+\tfrac{3}{2}{ R}_{111}) {
      R}_1 +(-\tfrac{9}{2} { R}_{222}
    \\
    & -\tfrac{27}{8} { R}_{122}+\tfrac{27}{8} { R}_{112}+\tfrac{3}{4}
    { R}_{111}) { R}_2+ \tfrac{15}{2}{ R}_{11}^2+(\tfrac{15}{2}{
      R}_{22}-\tfrac{45}{2} { R}_{12}) { R}_{11}+15 { R}_{12}^2-15 {
      R}_{12} { R}_{22})+\tfrac{1}{{ R}^2}((\tfrac{21}{8} { R}_{12}
    \\
    & -\tfrac{21}{8} { R}_{11}-\tfrac{21}{8} { R}_{22}) {
      R}_1^2+(\tfrac{21}{16} { R}_{22}-\tfrac{21}{16} {
      R}_{12}+\tfrac{21}{16} { R}_{11}) { R}_2 { R}_1+(\tfrac{63}{8} {
      R}_{22}-\tfrac{63}{8} { R}_{12}+\tfrac{63}{8} { R}_{11}) {
      R}_2^2)) s+(147 { R}_2+6 { R}_1){ R}^2
    \\
    & +(36 { R}_{122}-66 { R}_{112}+3 { R}_{111}) { R}+(30 {
      R}_{22}-27 { R}_{12}-39 { R}_{11}) { R}_1+(-45 {
      R}_{22}+\tfrac{33}{2} { R}_{12}-\tfrac{15}{4} { R}_{11}) {
      R}_2+9 { R}_{11122}
    \\
    & -6 { R}_{11222}-3 { R}_{11112}+\tfrac{1}{{ R}}(\tfrac{57}{8} {
      R}_1^3+\tfrac{3}{16} { R}_1^2 { R}_2+ (-3 { R}_{1122}+15 {
      R}_2^2+\tfrac{3}{2}{ R}_{1112}+\tfrac{3}{2}{ R}_{1222}) {
      R}_1+\tfrac{9}{4} { R}_2^3+(-\tfrac{21}{4} { R}_{1122}
    \\
    & +\tfrac{9}{2} { R}_{1222}+\tfrac{3}{4} { R}_{1112}) {
      R}_2+(-\tfrac{15}{2}{ R}_{122}+\tfrac{15}{2}{ R}_{112}) {
      R}_{11}+(-15 { R}_{112}+15 { R}_{122}) { R}_{12}) +\tfrac{1}{{
        R}^2}((\tfrac{21}{8} { R}_{122}-\tfrac{21}{8} { R}_{112}) {
      R}_1^2
    \\
    & +(\tfrac{21}{16} { R}_{112}-\tfrac{21}{16} { R}_{122}) { R}_2 {
      R}_1+(-\tfrac{63}{8} { R}_{122} +\tfrac{63}{8} { R}_{112}) {
      R}_2^2)
\end{alignat*}

\bigskip

\end{document}